\documentclass{article}

\usepackage{fullpage} 
\usepackage{amssymb} 
\usepackage{amsmath}
\usepackage{graphics}
\usepackage{epsfig}
\usepackage{pstricks}
\newtheorem{definition}{Def\hskip 1pt inition }[section]
\newtheorem{thm}[definition]{Theorem}

\newtheorem{defi}[definition]{Definition}
\newtheorem{prop}[definition]{Proposition}
\def\ds{\displaystyle}

\title{Comparison of numerical simulations of reactive transport\\ and chemostat-like models}
\author{{\sc I. Haidar${}^{1}$, F. G\'erard${}^{2}$ and A. Rapaport${}^{1}$\thanks{Corresponding author (Alain Rapaport)}}\\[2mm]
${}^{1}$ UMR INRA/SupAgro 'MISTEA' and EPI INRA/INRIA 'MODEMIC'\\
2, pl. Viala 34060 Montpellier, France\\
e-mails: {\tt haidar@supagro.inra.fr,rapaport@supagro.inra.fr}\\[2mm]
${}^{2}$ UMR INRA/SupAgro/CIRAD/IRD 'Eco\&Sols'\\
2, pl. Viala 34060 Montpellier, France\\
e-mail: {\tt gerard@supagro.inra.fr}
} 

\begin{document}

\maketitle

\noindent {\em Abstract.} 
The objective of the paper is to evaluate the ability of reactive transport models and their numerical implementations (such as MIN3P) to simulate simple microbial transformations in conditions of chemostat or gradostat models, that are popular in microbial ecology and waste treatment ecosystems. To make this comparison, we first consider an abstract ecosystem composed of a single limiting resource and a single microbial species that are carried by advection. In a second stage, we consider another microbial species in competition for the same limiting resource. Comparing the numerical solutions of the two models, we found that the numerical accuracy of simulations of advective transport models performed with MIN3P depends on the evolution of the concentrations of the microbial species: when the state of the system is close to a non-hyperbolic equilibrium, we observe a numerical inaccuracy that may be due to the discretization method used in numerical approximations of reactive transport equations.
Therefore, one has to be cautious about the predictions given by the models.

\section{Introduction}

The chemostat is a popular apparatus, invented simultaneously by Monod \cite{30} and Novick and Szilard \cite{33}, that allows the continuous culture of micro-organisms in a controlled medium. The chemostat has the advantage to study bacteria growth at steady state, in contrast to batch cultivation. 
The chemostat model serves also as a representation of aquatic natural ecosystems such as lakes. 
In the classical experiments involving chemostats, the medium is assumed to be perfectly mixed, that justifies mathematical models described by systems of ordinary differential equations \cite{14}. In natural ecosystems, ground-waters or industrial applications that use large reservoirs, the assumption of perfectly mixed medium is questionable, leading to spatialized models such as systems of nonlinear partial differential equations \cite{37}. However, nonlinear partial differential equations are difficult mathematical objects to understand, analyze and simulate. Even numerical schemes pose significant difficulties, particularly when solving coupled systems involving stiff reactions \cite{3, 34, 43}. Spatial considerations 
can be introduced in the ``classical'' model of the chemostat in simpler ways, as it is done in the gradostat model \cite{19} that mimics a series of interconnected chemostats of identical volumes. Other kinds of interconnection can be considered in order to cope with heterogeneity of porous media, considering stagnant water with small diffusion, mixing 
due to porosity \cite{13, 31}. If course, such representations are still oversimplified with regard to the complexity of natural porous media.

Over the past three decades, numerous reactive transport codes have been developed to study complex interactions between geochemical and transport processes in porous media. A number of reactive transport computer codes exist. Let us mention for instance the models COMEDIE-2D \cite{6}, CRUNCH \cite{38}, PHREEQC \cite{35}, ECOSAT \cite{16}, ORCHESTRA \cite{26} , RAFT \cite{4}, RT3D \cite{5}, HYTEC \cite{41}, HP1 \cite{15} and  MIN3P \cite{23}. To our knowledge, some of  these numerical tools, such as COMEDIE-2D and PHREEQC, are not suitable for unsaturated porous media and thus cannot be readily applied to soils (excepted in the special case of wetland soils). A range of other limitations can be found as well. For example, RT3D does not include equilibrium-controlled reactions, while ECOSAT neglects kinetically-controlled reactions and is limited to a single spatial dimension.

The standard approach for evaluating the computational accuracy of a reactive transport code is to compare its numerical results to those obtained from an analytical or a semi-analytical solution \cite{9, 34, 40, 42}. Unfortunately, analytical solutions are only available for simplified systems, such as the reactive transport of a single solute in 1-D homogeneous systems at steady state, which is well behind the actual capabilities of the models. To remediate to this, inter-comparison of numerical codes has been largely employed. This inter-comparison involves the independent solution of the same problem using a variety of numerical techniques \cite{3, 7, 10, 25, 6}.

This study aims at comparing the accuracy of a reactive transport model with other kinds of models such as the mathematical model of the chemostat. This confrontation takes place in the framework of microbial ecology, for which concepts of competition and coexistence are crucial \cite{1, 2, 8, 36, 39, 39b}. We have chosen the reactive transport code MIN3P \cite{23} for this study. This reactive transport model is notably recognized for its numerical robustness \cite{3,25}. In addition, the model MIN3P can simulate general reactive transport problems in variably saturated media for 1D to 3D systems. The flow solution is based on Richard's equation \cite{32}, and solute transport is simulated by means of the advection-dispersion equation. Gas transport can be taken into account as well, either by considering advection and Fick diffusion or the Dusty Gas Model \cite{27, 28}. A range of bio-geochemical processes are included in MIN3P (aqueous speciation, mineral dissolution-precipitation, gas dissolution/exsolution, ion exchange, and competitive or non competitive sorption). A generalized kinetic expression for dissolution-precipitation and intra-aqueous reactions allows the consideration of fractional order or Monod-type rate expressions, and parallel reaction pathways. This code has been used for a number of applications in different fields of environmental science, ranging from inorganic and organic contaminant transport and groundwater remediation \cite{21, 22, 24, 29} to soil hydrology and bio-geochemical cycles in terrestrial ecosystems \cite{11, 12, 20}.

The chemostat model is the simplest mathematical model for describing the dynamics of microbial growth under a constant flow of substrate, and its theory is well understood \cite{37}. In this work, we consider a series of interconnected chemostats for the simulation a one-dimensional heterogeneity, that we compare with the solutions provided by reactive transport models considering the same spatial structure.

\section{Material and method}

\begin{figure}[!ht]
\begin{center}
\begin{pspicture}(1,0)(5,3)
\put(-0.8,0){\includegraphics[width=8cm,height=3cm,angle=0]{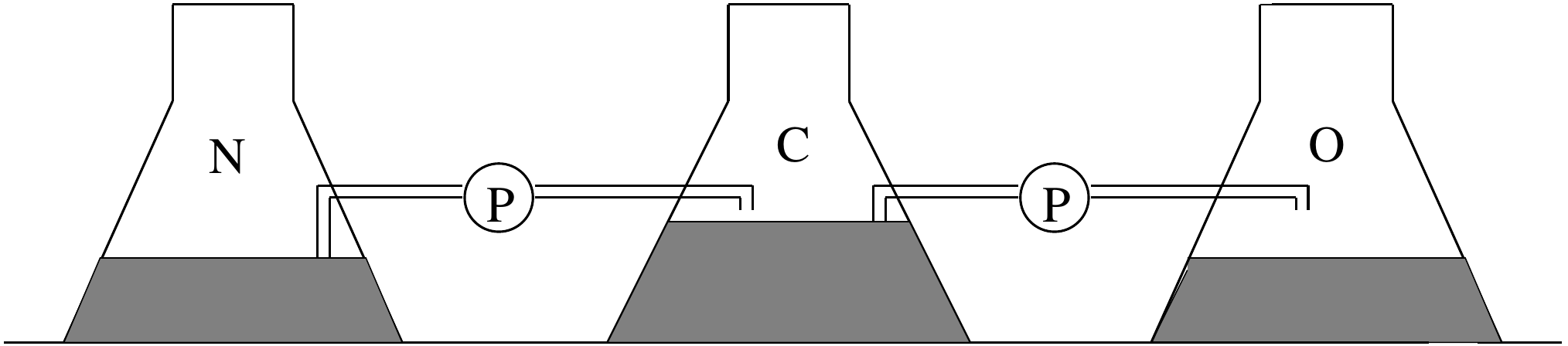}}
\end{pspicture}
\caption{A schematic view of the chemostat experiment}
\label{labchem}
\end{center}
\end{figure}

In practice, a chemostat in laboratory is an apparatus that consists of three connected vessels as shown in Fig. \ref{labchem}. The leftmost vessel is called the feed bottle and contains all of the nutrients needed for the growth of a microorganism. The central vessel is called the culture vessel, while the last is the overflow or collection vessel. The content of the feed bottle is pumped at a constant rate into the culture vessel, while the content of the culture vessel is pumped at the same constant rate into the collection vessel. We denote by $S_{in}[mol/l ]$ the constant concentration of nutrient pumped with a volumetric flow rate that we denote by $Q [l.s^{-1}]$. The dilution rate $[s^{-1}]$ is defined as $D=Q/V,$ where $V$ is the volume of the culture vessel. We shall also denote by $\mu(.)[s^{-1}]$ the specific growth rate of micro-organisms and $k$ the yield factor of the bio-conversion. The dynamical model of the chemostat can then been written in the following way.

\begin{equation}\label{chemostat}
\left\{\begin{array}{lll}
 &\displaystyle{\frac{dS}{dt}=-\frac{\mu(S)}{k}B+D(S_{in}-S)}\vspace{5 pt} \cr&
\displaystyle{\frac{dB}{dt}=(\mu(S)-D)B}\vspace{5 pt}
\end{array} \right.
\end{equation}   
                                                                  
Because of the boundary conditions (i.e. input of nutrients in the culture vessel and output of contents from the culture vessel), numerical implementations of reactive transport models such as MIN3P are not able to simulate straightforwardly a single ecosystem such as the one in the culture vessel. Indeed, the use of the logarithm of the concentrations in the numerical code prevents to having a null concentration of biomass at the input of the culture vessel.
In order to simulate an ecosystem in a single tank, one has to consider three control volumes. 
 Moreover, for intrinsic reasons, three control volumes is the minimum number for a one dimensional discretization in numerical implementations 
of reactive transport models, such as the MIN3P code that we have chosen to simulate our ecosystem.
In this way, the numerical implementation is closer from the true laboratory experiment with three vessels, that we described above.
Nevertheless, we shall refer in the following to the {\em chemostat} for the culture vessel only.

With MIN3P we began by simulating a simple example of three chemostats in series having the same volume, in presence of a single biomass $B$ and a single substrate $S$ with the same initial conditions in the three control volumes. The dynamical model representing $n$ chemostats connected in series is given by the equations:                     
\begin{equation*}\label{system}
\left\{\begin{array}{lll}
 &\displaystyle{\frac{dS_i}{dt}=-\frac{\mu(S_i)}{k}B_i+D_i(S_{i-1}-S_i)}\vspace{5 pt} \cr&
\displaystyle{\frac{dB_i}{dt}=(\mu(S_i)-D_i)B_i}+D_iB_{i-1},\cr 
\end{array} \right.
\end{equation*}                                   
where $D_i= Q/V_i, S_i (resp B_i)$ represents the substrate concentration (resp biomass concentration) in the $ith$ bioreactor $(i=1,…,n).$ $S_0=S_{in}$ and $B_0=0$ and $V_i=\frac{V}{n}$.\\
We consider that the qualitative behavior of this system of ordinary differential equations (ODE) is today well understood (see the Appendix), and we used different robust numerical schemes for solving this system of ODE (Runge-Kutta, LSODA ...) that were all thoroughly consistent with the analytic analysis of the steady states and their stability. That amounts to assume that we can trust the numerical solutions obtained by the numerical integration of ODEs for this model.\\
 
\begin{figure}[!ht]
\begin{center}
\begin{pspicture}(-4,0)(0,7.7)
\put(-5,0.5){\includegraphics[width=5cm,height=7cm,angle=0]{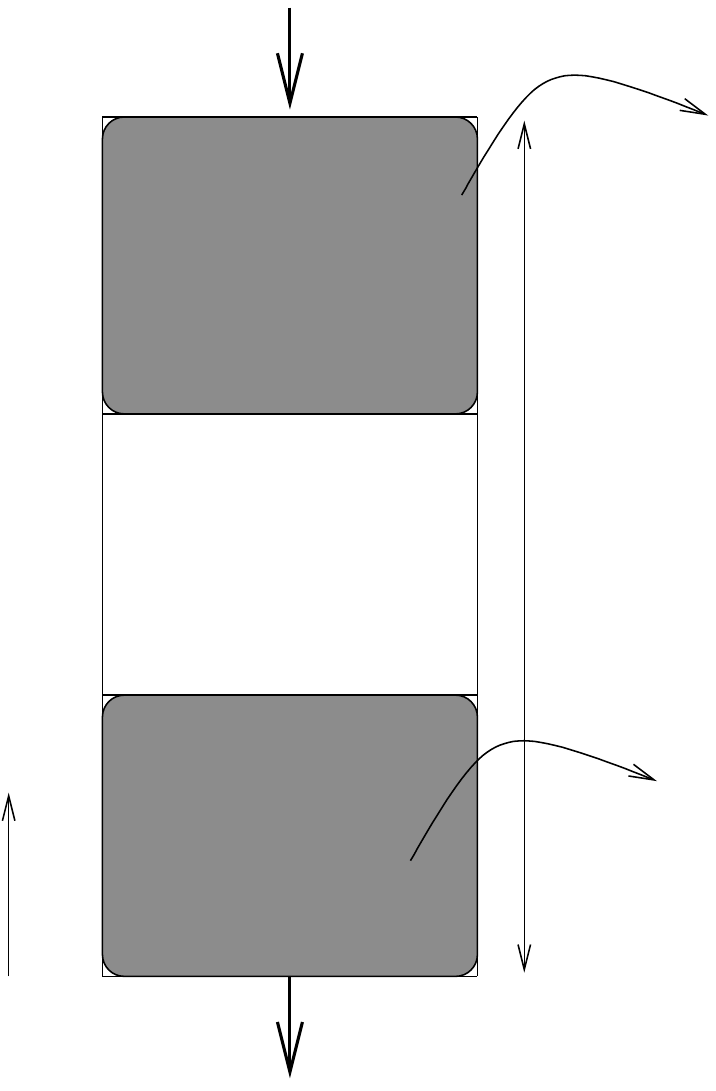}}
   \rput(-2.9,0.2){$S_{out}$}
   \rput(-3.3,7){$Q$}
   \rput(-3,7.7){$S_{in}$}
   \rput(-3.3,0.9){$Q$}
   \rput(0.5,6.5){Third C.V.}
   \rput(-0.3,4){$1m$, 3 cells}
   \rput(0.5,2.2){First C.V.}
   \rput(-5.3,1.5){Z}
\end{pspicture}
\caption{Configuration of three control volumes}
\label{}
\end{center}
\end{figure}
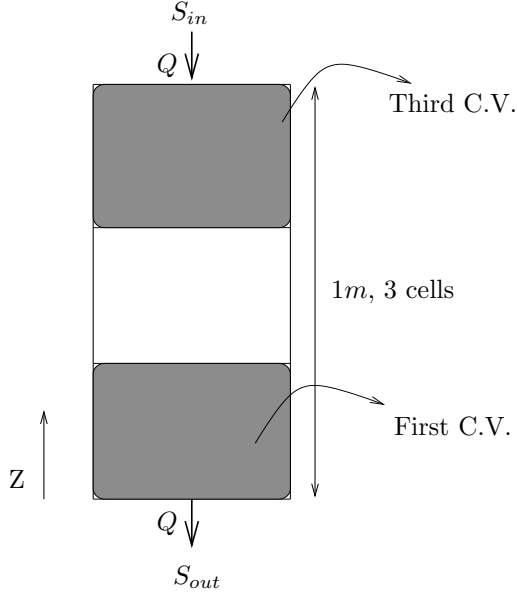

For the simulation of the chemostat model with reactive transport models, we consider a boundary domain in three dimensions. The boundary conditions for the liquid flow are of first type with a value of $0$ in the output face, and of second type specifying a flux of $Q [l.s^{-1} ]$ in the input and output faces. A specific choice of flow condition gave us a fully water saturated medium at any time . The boundary conditions for the reactive part are of second type in the output face and of third type in the input face with a value of the substrate concentration equal to $S_{in}$. 

To simulate several chemostats connected in series without diffusion, we have chosen the free diffusion coefficient in water and air, as well as the specific storage coefficient equal to zero. The porosity of the medium formed by one domain only is chosen equal to one. The day has been chosen as the time unit, with a maximum time increment of $10^{-3}$ day and a starting time step of $10^{-10}$ day. 

Finally, we have introduced a simple theoretical irreversible reaction expressing the bio-conversion of moles of substrate into one mole of biomass. The specific growth rate of biomass considered here is of Monod's type:

$$\mu (S) = \frac{\mu_{max}S}{k_s + S},$$  where $\mu_{max}[s^{-1}]$ represents the maximum of the intra-aqueous kinetic reaction and $k_s[mol/l]$ the half saturation constant.\\

In a second stage, we consider the same spatial considerations but with two species instead of one, assuming that each species has a growth function that follows a Monod law, as described above. 
We assume that their interaction is due only to a common limiting resource, that is the species compete for the same substrate. We focus on the case of a ``true'' competition, in the sense that we assume 
that a species is the most efficient one when the resource is very rare, while the other species is better when the resource is less rare. 

\section{Results and discussion}

Denote by $S_M ^*$ (resp. $S_C ^*$) the value of the substrate concentration computed by MIN3P (resp. by the chemostat model) at the steady state, and by $B_M ^*$ the value of the biomass concentration at steady state given by MIN3P. Let $\delta$
be the absolute value of the difference between $S_M ^*$ and $S_C ^*$, that serves as an indicator of divergence between simulations of both models.\\

For the comparison of the two models, we first represent the indicator $\delta$ as a function of $Q$. For the study of the effect of a spatial discretization, we shall then consider $\delta$ as a function of the number of cells. \\

For convenience, we shall omit the units of the numerical values that we give below (concentrations are assumed to be measured in $mol.l^{-1}$, volumes in $l$, flux in $l.s^{-1}$ and growth rates in $s^{-1}$.

\subsection{The single species case}

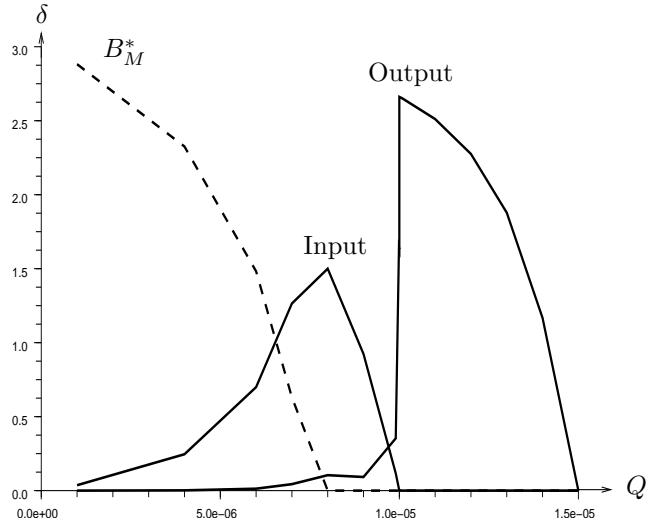
\begin{figure}[!ht]
\begin{center}
\begin{pspicture}(0,0)(2,7.5)
\put(-2.8,0.5){\includegraphics[width=8cm,height=6.5cm,angle=0]{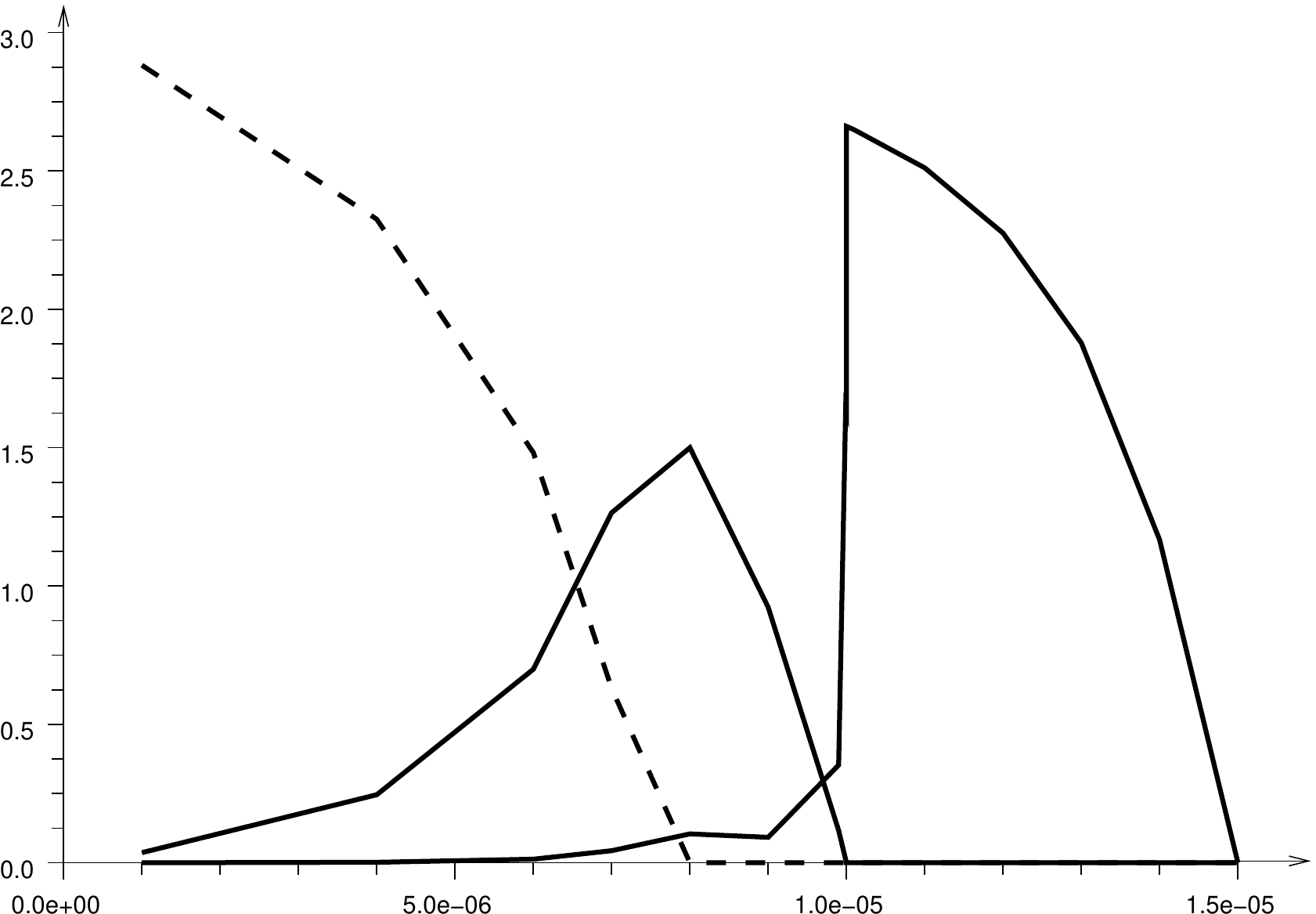}}
   \rput(5.5,0.9){$Q$}
   \rput(-2.4,7.2){$\delta$}
   \rput(-1.3,6.7){$B_{M}^{*}$}
  \rput(1.5,4.1){Input}
 \rput(2.5,6.4){Output}
\end{pspicture}
\caption{Comparison for three control volumes}
\label{comp3}
\end{center}
\end{figure}

Choosing $\mu_{max}=4.10^{-5}$, we have studied the variation of $\delta$ with respect to $Q$ in the first and third control volumes. Each  volume has been chosen equal to $1$. On the graphs of $\delta$ that correspond to the first (resp. third) control volume (see Figure \ref{comp3}, Output (resp. Input)), we notice that for 
$10^{-6}\leq Q \leq10^{-5}$, one faces a significant difference between simulations of both models in the input control volume, and that is even greater in the output control volume. But for $10^{-5}\leq Q \leq15.10^{-6}$, one has almost no difference in the input control volume, and one can observe a jump of $\delta$ about $Q=10^{-5}$ in the output control volume. For the simulations, we have chosen $S_{in}=3,k_s=1,B(0)=2$ and $S(0)=5.$\\

We know from the theory of the chemostat (see Appendix, Proposition
\ref{prop1species}), that for a "perfectly mixed" tank, with a single species growing on a single limiting substrate, the condition $\frac{\mu(S_{in})}{D}>1$ ensures that the biomass $B$ is not wash-out. This result is surprisingly not observed in simulations with MIN3P. To show that, we studied the variation of the biomass concentration $B$ with respect to $Q$ in the input control volume. We notice that the biomass is washed out for $Q\geq 8.10^{-6}$ (see Figure \ref{comp3}, $B_M^*$ in dashed line). But when $8.10^{-6}\leq Q<10^{-5}$, we have $\frac{\mu(S_{in})}{D_1} =\frac{10^{-5}}{Q}>1$, 
and one can observe on the graph of $\delta$ and the difference between MIN3P and the chemostat simulations in the input control volume. 

For the case of three chemostats connected in series with the same volume and traversed by the same flow rate $Q$, the removal of the biomass in the input control volume has to lead theoretically to its removal in the output control volume. But this is not the case with MIN3P and we can notice on Fig. \ref{comp3} (Output), that for 
$8.10^{-6}\leq Q \leq 15.10^{-6}$, we obtain the wash-out of biomass in the input control volume, the biomass in the output control volume being not yet washed out. In other words, under certain conditions, the microbial growths predicted by both simulations are radically different.\\

\begin{figure}[!ht]
\begin{center}
\begin{pspicture}(0,0)(2,7.5)
\put(-2.8,0.5){\includegraphics[width=8cm,height=6.5cm,angle=0]{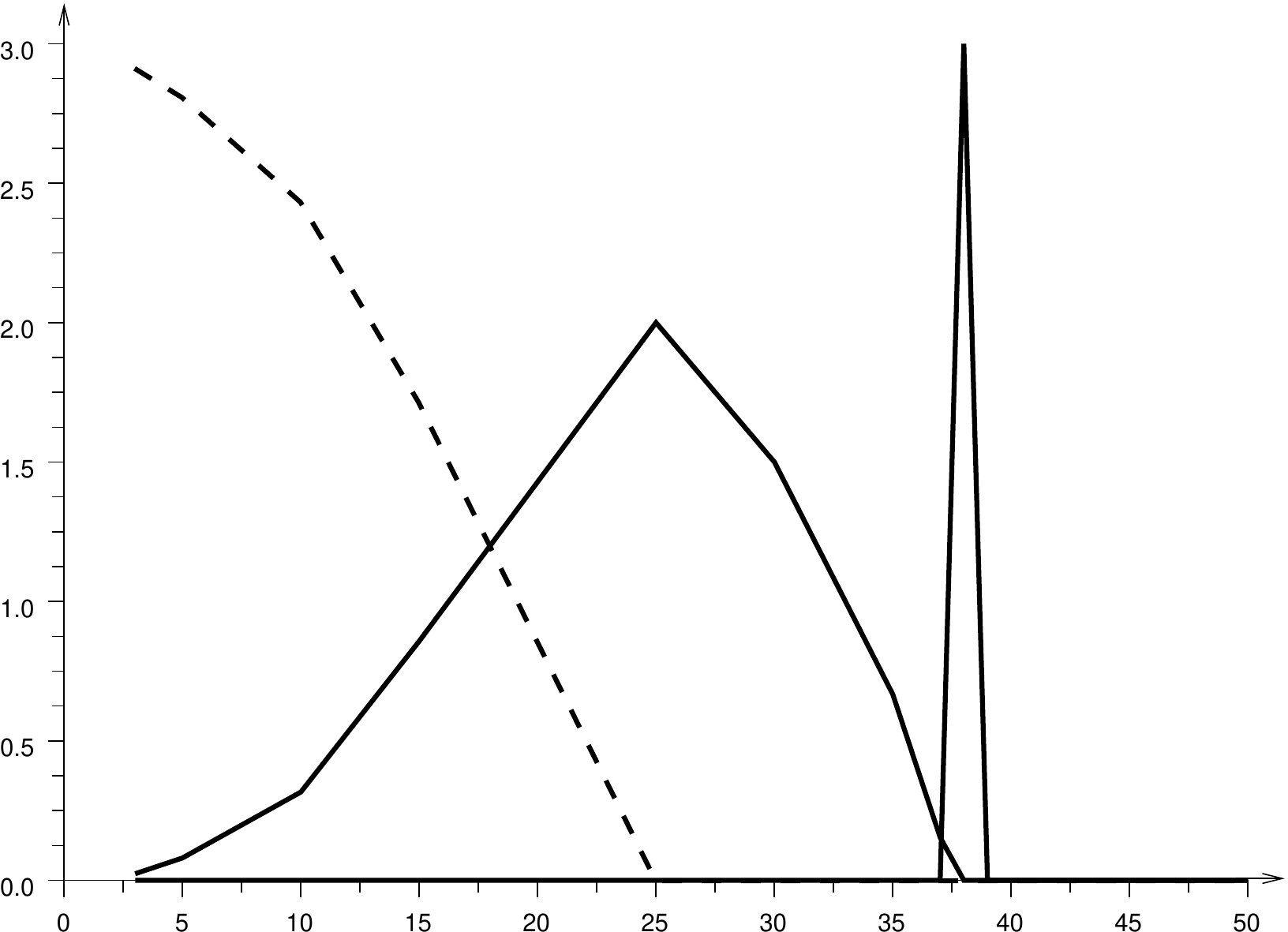}}
   \rput(0.3,0.2){Number of control volume}
\rput(-2.4,7.2){$\delta$}
   \rput(-1.3,6.7){$B_{M}^{*}$}
  \rput(1.1,5.1){Input}
 \rput(3.1,7){Output}
\end{pspicture}
\caption{Comparison with several control volumes}
\label{compn}
\end{center}
\end{figure}
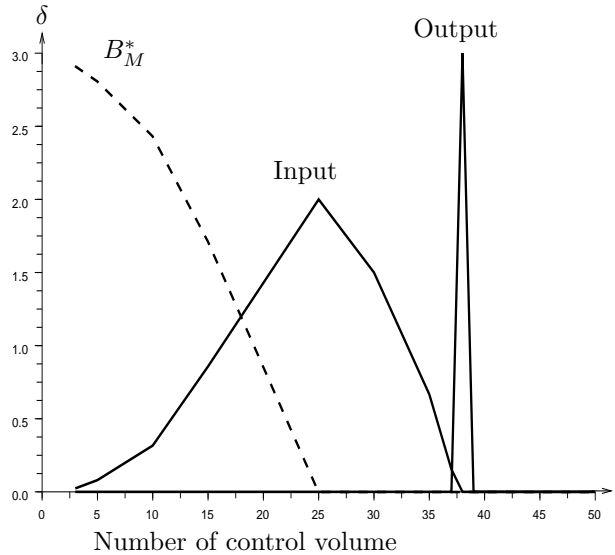

For the study of the effect of a spatial discretization on the difference between numerical reactive transport and chemostat models at steady state, we took the same conditions as before, with a maximal kinetic rate equal to $\mu_{max}=5.10^{-4}$ and a flow $Q = 10^{-5}.$ We did vary the number of discretization steps between $n=3$ and $n=50$ and studied the variation of $\delta$ with respect to $n$ (the total number of cells) in the input and output control volumes. 
Denote by $D_{n}$ the dilution rate in each control volume for a discretization in $n$ cells.
On the graphs of $\delta$ that corresponds to the input (resp. output) control volume (see Fig. \ref{compn}, Input (resp. Output)), we notice for $3\leq n \leq37$ (that theoretically corresponds to the survival of biomass, because one has $\frac{\mu(S_{in})}{D_n}=\frac{75}{2n}>1$) 
that we have the same behavior of $\delta$ as before. Similarly, for $39\leq n\leq50$ (values that correspond to  $\frac{\mu(S_{in})}{D_n}<1$), we have no difference between MIN3P and the chemostat model in the input and output control volumes. But for $n=38$, we observe the same jump of $\delta$ in the output control volume as previously observed.\\

We have observed on this example a significant difference between reactive transport and chemostat models when passing from a steady state of survival of the biomass to the wash out steady state. 
this corresponds to a bifurcation passing from two equilibriums (wash-out and biomass survival) to a single equilibrium (wash-out). The limiting case corresponds to a non-hyperbolic equilibrium (see the Appendix for the definition of hyperbolic equilibriums). Because of the use of the logarithm in the MIN3P code, we expect the internal solution to take very large values when the concentration of biomass tends to zero. One can also 
detect this phenomenon on the simulations when noticing a "time lag" between MIN3P and chemostat trajectories about the steady state.

\subsection{The two species case}

To emphasize the problem that occurs about non-hyperbolic equilibrium in the case of one species, we present in this part a more subtle 
situation, considering two species that compete for the same substrate. The extension of the model (\ref{chemostat}) is given by the following equations:
\begin{equation}\label{comp}
\left\{\begin{array}{lll}
 \ds \frac{dS}{dt}&=&\ds -\frac{\mu_1(S)}{k_1}B_1-\frac{\mu_2(S)}{k_2}B_2+D(S_{in}-S)\vspace{5 pt} \cr
\ds \frac{dB_1}{dt}&=&(\mu_1(S)-D)B_1\vspace{5 pt} \cr
\ds \frac{dB_2}{dt}&=&(\mu_2(S)-D)B_2\vspace{5 pt} \cr 
\end{array}\right. 
\end{equation}
where functions $\mu_1(.)$ and $\mu_2(.)$ denote the kinetics of species $B_1$ and $B_2$ respectively. 
For this system, the
Proposition \ref{prop2species} given in the Appendix shows that 
non-hyperbolic equilibrium points could exists away from the wash-out equilibrium. For the one species case, a non-hyperbolic 
equilibrium could exists but on the boundary of the positive domain only. 
For the two species case, many non-hyperbolic equilibriums could exists on the interior of the positive domain.\\

Let $\lambda_i$ be the positive solution, when it exists, of $\mu_i (S)=D$ for $i=1,2$. Under the condition $S_{in}\geq \max(\lambda_1,\lambda_2)$, 
the system (\ref{comp}) admits generically three equilibrium points, given by $E_0 (S_{in},0,0), E_1 (\lambda_1,k_1 (S_{in}-\lambda_1),0)$ and $E_2 (\lambda_2,0,k_2 (S_{in}-\lambda_2))$.
We distinguish now two different cases. 
If for some $i=1,2$ one has $\lambda_i=S_{in}$, then the
equilibrium $E_0$ is non-hyperbolic as before (see Proposition \ref{prop1species} in the Appendix). Moreover if one has $\lambda_1=\lambda_2<S_{in}$, then $E_1=E_2$ is also a non hyperbolic equilibrium (Proposition \ref{prop2species} in the Appendix). Starting from our observations in the case of one species, we have built suitable examples to study the behavior of MIN3P about those non-hyperbolic equilibriums. \\

To understand the comparison, we first recall that the mathematical theory of the model (\ref{comp}) predicts the competitive exclusion in the generic case, that is at most one competitor avoids the extinction (see the Appendix).
This property refers to the well-known Competitive Exclusion Principle in ecology, that has been widely studied in the literature
(see for instance \cite{1, 2, 8, 39}).
For the chemostat model, the Principle can be stated as follows.

Considering two increasing growth rates $\mu_1 (.)$ and $\mu_2 (.)$ such that both $\lambda_1$ and $\lambda_2$ are smaller than  $S_{in}$ (for a sufficiently large $S_{in}$). Then, one has the following issue of the competition for large times
\begin{itemize}
\item[-] when $\lambda_1<\lambda_2$, the species $B_1$  avoids the extinction, 
\item[-] when $\lambda_1>\lambda_2$, the species $B_2$  avoids the extinction. 
\end{itemize}
For the non generic case $\lambda_1=\lambda_2$, 
it is possible to predict the coexistence of the two species, 
invalidating the Principle (on this single chemostat case).\\

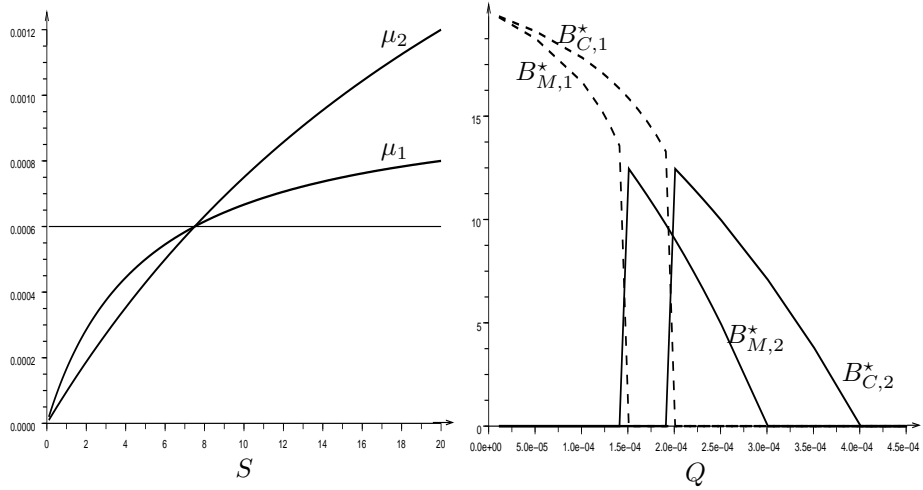
\begin{figure}[!ht]
 \begin{center}
 \begin{pspicture}(0,0)(7,7)
\put(3.5,7){\includegraphics[width=6cm,height=6cm,angle=-90]{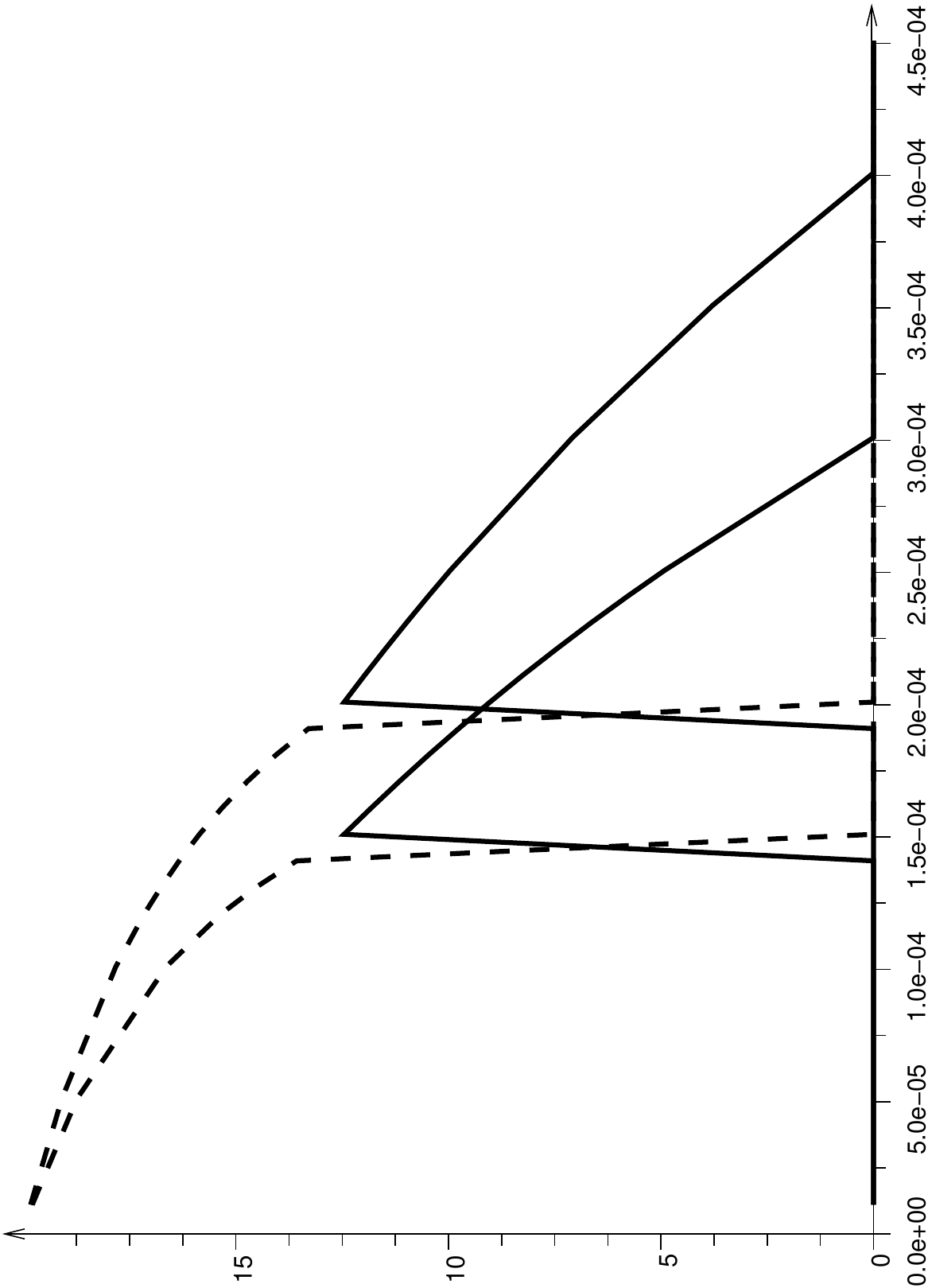}}
\put(-3,7.5){\includegraphics[width=7cm,height=7cm,angle=-90]{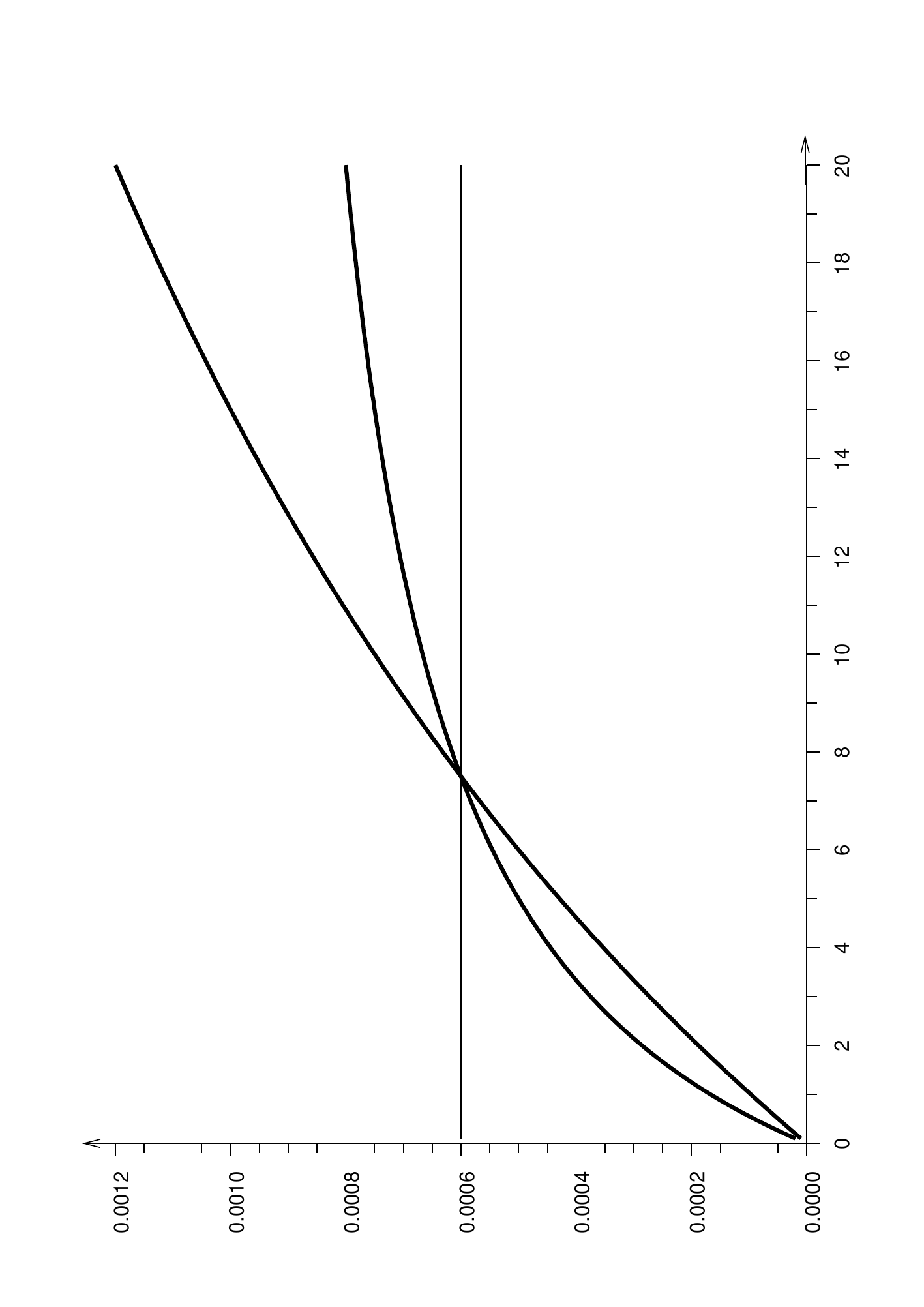}}
\rput(4.5,6){$B^\star_{M,1}$}
\rput(7.3,2.5){$B^\star_{M,2}$}
\rput(5,6.5){$B^\star_{C,1}$}
\rput(8.8,2){$B^\star_{C,2}$}
\rput(2.5,5){$\mu_1$}
\rput(2.5,6.5){$\mu_2$}
\rput(6.5,0.7){$Q$}
\rput(0.5,0.8){$S$}
 \end{pspicture}
\caption{Two species in competition in a single chemostat}
\label{comp1}
 \end{center}
 \end{figure}

In the simulations, we have considered two species $B_1$ and $B_2$ 
with a specific growth rate given by
$$\mu_1(S)=1.10^{-3}\frac{S}{5+S} \quad \mbox{and} \quad \mu_2(S)=3.10^{-3}\frac{S}{30+S} \ . $$
One can notice on Fig. \ref{comp1} (left) that the graphs of these two functions intersect away from zero. This implies that depending on the dilution rate $D$, the corresponding value of $\lambda_{1}$ can be less or greater than $\lambda_{2}$.
The input concentration $S_{in}$ has been chosen equal to $20$ and the initial state vector has been kept equal to
$(S(0), B_1 (0), B_2 (0))=(5,2,3)$. A simple calculation shows that for $Q=2.10^{-4}$, one has $\lambda_1=\lambda_2=\frac{15}{2}<S_{in}$. Then for this choice of $\mu_1, \mu_2$ and $Q$ the dynamical system (\ref{comp}) admits positive non-hyperbolic equilibriums. We aim to study the numerical evolution of the species concentrations about those equilibriums in both models, depending on the choice of the flow $Q$. For this purpose, we have plotted the graphs of the species concentrations at steady state given by MIN3P and the chemostat model in the input control volume (see Fig. \ref{comp1}, right), denoting
by $B_{M,i}^*$ (resp. $B_{C,i}^*)$ for species $i=1,2$. We observe that for $Q\leq10^{-5}$, both simulations present almost the same solutions. When $Q$ increases a difference between the models begins to appear until we detect a wrong prediction of the species that avoids the extinction. This happens for $15.10^{-5}\leq Q\leq 4.10^{-4}$. For $Q\geq4.10^{-4}$, both simulations show again almost the same solutions. \\

So we have observed another significant difference between simulations of reactive transport and chemostat models about a bifurcation point, where the species that avoid extinction switches. \\

To study the effect of a spatial discretization in the case of two species, we have chosen a specific example of twenty perfectly-mixed tanks that are interconnected in series. The volumes $V_i (i=1, ..., 20)$, the input concentration $S_{in}$, the flux $Q$ and the specific growth rates $\mu_1(.)$ and $\mu_2(.)$ are chosen as follows, in such a way that the species $B_1$ passes from the wash-out state (in $V_1,V_2 \mbox{ and } V_3)$ to a coexistence state (in $V_i, i=4, ...,20$). 
\vspace{10pt}

$
\left.\begin{array}{lll}
V_1=10.10^{-2}\vspace{5 pt} \cr
V_2=9.10^{-2}\vspace{5 pt} \cr
V_3=1.10^{-2} \vspace{5 pt}\cr 
V_4=11.10^{-2}\vspace{5 pt} \cr
V_5=9.10^{-2} \vspace{5 pt}\cr 
V_6=\cdots=V_{20}=4.10^{-2}\cr
\end{array} \right.
$
\hspace{10pt}
$
\left.\begin{array}{lll}
\ds \mu_1(S)=4.629.10^{-5}\frac{S}{6+S}\vspace{5 pt} \cr
\ds \mu_2(S)=6.944.10^{-5}\frac{S}{18+S}\vspace{5 pt}\cr 
 Q=0.3587.10^{-5}\vspace{5 pt} \cr
 S_{in}=19.25\vspace{5 pt} \cr
%t_{min}=10^{-8}, t_{max}=10^{-4} \cr
\end{array} \right.
$\\

Under MIN3P we have discretized the domain into twenty control volumes, and have compared the solutions of both models in each control volume at steady state.

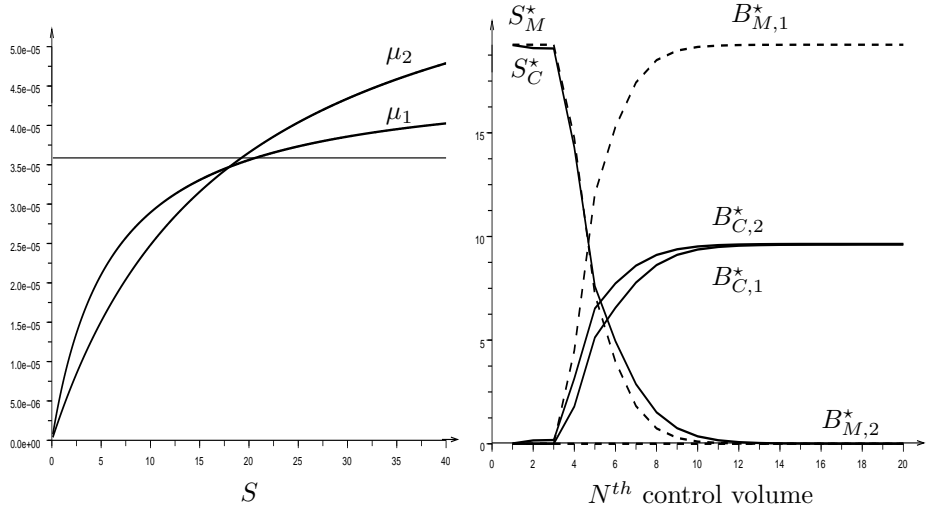
\begin{figure}[!ht]
 \begin{center}
 \begin{pspicture}(0,0)(7,7)
\put(3.5,6.5){\includegraphics[width=6cm,height=6cm,angle=-90]{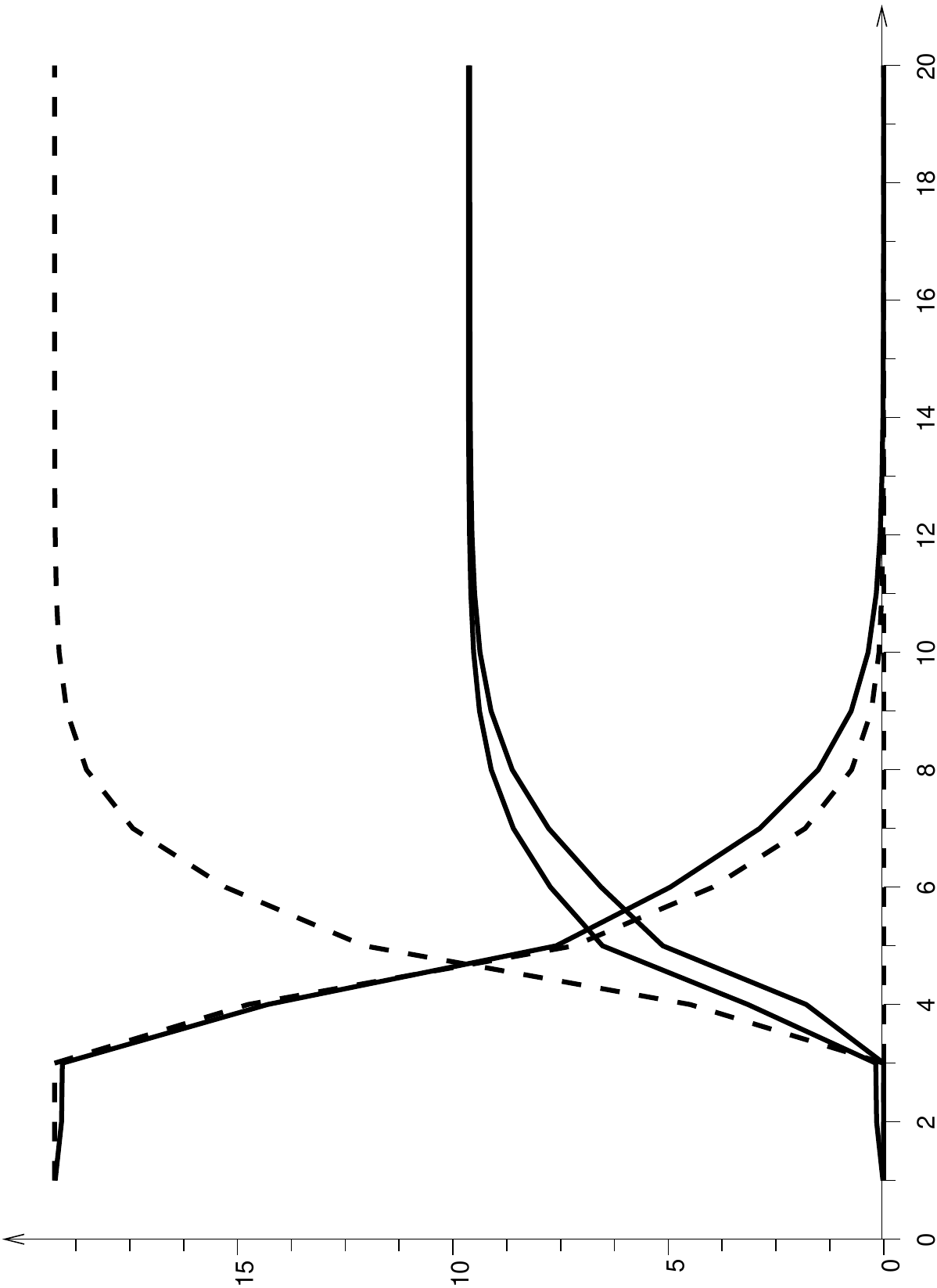}}
\put(-3,7){\includegraphics[width=7cm,height=7cm,angle=-90]{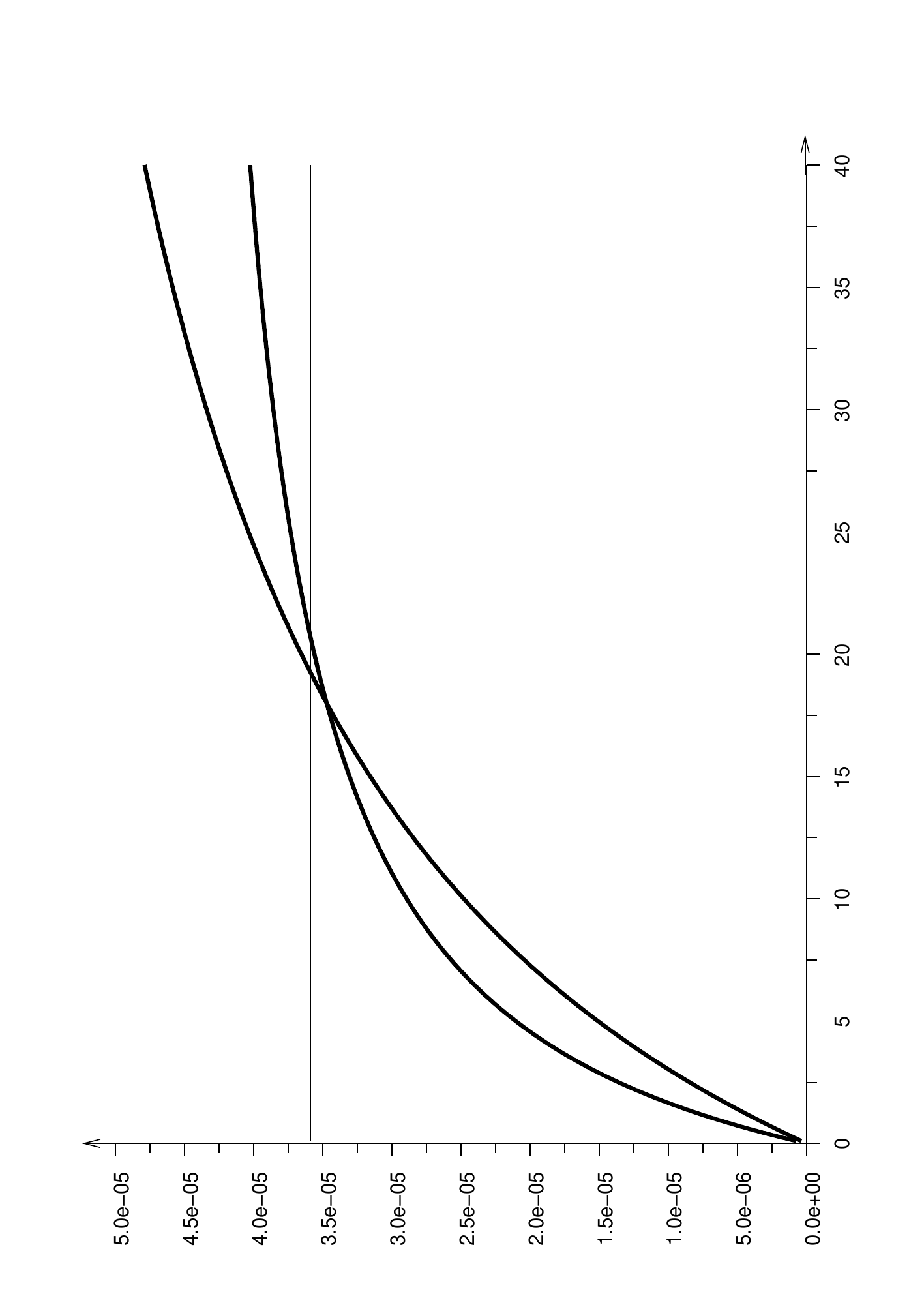}}
\rput(2.5,5.2){$\mu_1$}
\rput(2.5,6){$\mu_2$}
\rput(4.2,6.5){$S^\star_{M}$}
\rput(7.3,6.5){$B^\star_{M,1}$}
\rput(8.5,1.1){$B^\star_{M,2}$}
\rput(4.2,5.8){$S^\star_{C}$}
\rput(7,3){$B^\star_{C,1}$}
\rput(7,3.8){$B^\star_{C,2}$}
\rput(6.5,0.2){$N^{th}$ control volume}
\rput(0.5,0.2){$S$}
 \end{pspicture}
\caption{Two species in competition in a series of chemostats}
\label{comp2}
 \end{center}
 \end{figure}
On Figure \ref{comp2} (right), one can observe that the substrate concentrations computed at steady state with both models are almost
identical in the first control volumes, up to the third one 
where the solutions are different, and become radically divergent in further control volumes. Here, the solution computed by MIN3P does not predict the coexistence of two species...
Let us underline that the species coexistence is no longer a pathological case when chemostats of different volumes are connected in series (see for instance \cite{39b,36}).\\

To summarize, we have shown that under certain circumstances the issue of the winner of the competition between the two species are 
predicted radically different by the simulations of both models.

\section{Conclusion}
Our results show a possible inaccuracy of numerical reactive transport models 
in the simulation of the dynamics of simple ecosystem of chemostat type. Our objective was not to challenge the model MIN3P for its ability to simulate complex problems involving mass flow and multicomponents reaction networks, but we raise the fact that the numerical accuracy of the model MIN3P depended on the evolution of microbial species. When the hydrodynamic conditions make the system close to a washing-out of one or several microbial species, or to the coexistence of species, we observe numerical bias in the computation of the solution that leads to radically different predictions. Consequently one may wonder if the numerical issue found here in simple systems (a single substrate and one or two microbial species, advective transport) prevails in more complex systems, when a network of kinetically-controlled reactions is considered to simulate for instance remediation problems in ground-waters. We believe that a study of the eigenvalues of the linearized dynamics about steady states is important, for detecting possible numerical inaccuracy.

\section*{Appendix}

For the study of the behavior of autonomous dynamical systems in $ \mathbf{R}^n,$
\begin{equation}\label{nm}
\frac{dX}{dt}=F(X),
\end{equation}
with $F\in C^1(\mathbf{R}^n)$, one usually determine first its equilibrium points (denoted $X^{\star}$) as solutions of $f(X)=0$,
and then study the eigenvalues of the linear dynamics
\[ \frac{dX}{dt}=J(X^{\star})X,\]
called the linearization of (\ref{nm}) about $X^{\star}$, where $J(X^{\star})$ is the Jacobian matrix of $F$ at $X^{\star}$. If all the eigenvalues $v_i$ ($i=1,... n$ ) of $J(X^{\star})$ have nonzero real parts, then we said that $X^{\star}$ is hyperbolic. When at least one of its 
eigenvalues has a zero real part, then we said that $X^{\star}$ is non-hyperbolic (see for instance the reference \cite{17}).\\

We recall now the usual definitions of stability and a main result 
allowing to conclude about the nature of an equilibrium (see the
reference \cite{17} for more details).\\ 

\begin{defi}
The equilibrium $X^{\star}$ is said to be stable if for every $\epsilon>0$ there exists $\delta>0$ such that
$$\left\|X(t_0 )-X^{\star}\right\|<\delta\Rightarrow \left\|X(t)-X^{\star}\right\|<\epsilon \hspace{5pt} \forall t\geq t_0.$$
If this condition its not satisfied, the equilibrium is unstable.
\end{defi}

\begin{defi}
The equilibrium $X^{\star}$ is said to be (locally) exponentially stable if for every $\epsilon>0$ there exists three real numbers $a>0,b>0$ and $\delta>0$ such that
$$\left\|X(t_0 )-X^{\star}\right\|<\delta\Rightarrow \left\|X(t)-X^{\star}\right\|<a\left\|X(t_0)-X^{\star}\right\|e^{-bt} \hspace{5pt} \forall t\geq t_0.$$
\end{defi}

\begin{thm}
$\,$
\begin{itemize}
\item If all the eigenvalues of $J(X^{\star} )$ have negative real parts then the equilibrium point $X^{\star}$ is exponentially stable.
\item If one of the eigenvalues of $J(X^{\star})$ has a positive real part then the equilibrium point $X^{\star}$ is unstable.
\end{itemize}
\end{thm}

\noindent {\em Remark.} If the Jacobian matrix $J(X^{\star})$ has at least one eigenvalue with zero real part, then we need to use other results to conclude about the behavior of the trajectory of the system (\ref{nm}).\\

\medskip	

For the chemostat model (\ref{chemostat}), one has the following result.	
\begin{prop}
\label{prop1species}
Denote by $\lambda$ the solution when it exists of $\mu(S)=D$. 
\begin{itemize}
\item If $D<\mu(S_{in})$, the system (\ref{chemostat}) admits two equilibrium points given by $E_1(S_{in},0)$, which is unstable and $E_2(\lambda,k(S_{in}-\lambda))$, which is locally exponentially stable. 
\item If $D\geq \mu(S_{in})$, the only non-negative equilibrium point it is $E_1$ which is locally exponentially stable, excepted for the case $D=\mu(S_{in})$ for which it is non-hyperbolic.
\end{itemize}
\end{prop}

\textit{Proof.} We can easily verify that for any $t\geq0$, the trajectories of (\ref{chemostat}) remains in the first positive orthant,
and are bounded: one can straightforwardly write
\[
\frac{dB}{dt} +k\frac{dS}{dt}=D(kS_{in}-B-kS)
\]
from which one deduces that $t\mapsto B(t)+kS(t)$ is bounded that 
the trajectories of the system are bounded.
Determining the equilibrium points of the system (\ref{chemostat})
amounts to solve the following system  
\begin{equation}
\left\{\begin{array}{lll}
\ds -\frac{\mu(S)}{k}B+D(S_{in}-S)&=&0\vspace{5 pt} \cr
(\mu(S)-D)B&=&0.\cr
\end{array} \right.
\end{equation} 
The wash-out equilibrium point  $E_1 (S_{in},0)$ is always solution,
and there is a possibility of another equilibrium point
$E_2 (S^{\star},k(S_{in}-S^{\star}))$ with 
$S^{\star}=\lambda$, when $S^{\star}<S_{in}$.
To study the stability of these equilibrium points,
we write the Jacobian matrix of the system (\ref{chemostat}):

$$J(S,B)=\begin{pmatrix} -\dfrac{\mu'(S)}{k}B-D & -\dfrac{\mu(S)}{k}\vspace{7pt} \\
                          \mu'(S)B   & \mu(S)-D &
         \end{pmatrix}, $$

whose eigenvalues are $$v_1(S,B)=-D<0$$ and $$v_2(S,B)=\mu(S)-D-\dfrac{\mu'(S)}{k}B.$$
At $E_1$, one has $v_2(E_1)=\mu(S_{in})-D$ and at the non-trivial equilibrium $E_2$(when it exists), one has 
$v_2 (E_2)=-\frac{\mu^{'}(S^{\star})}{k} B^{\star}<0.$ 

So, when $\mu(S_{in})<D$, we conclude that the non-trivial equilibrium does not exist and we obtain that $E_1$ is locally exponentially stable. When $\mu(S_{in})>D$, $E_1$ is unstable and $E_2$ is locally exponentially stable. 
For the particular case $\mu(S_{in})=D$, the non-trivial equilibrium does not exist and $E_1$ is a non-hyperbolic equilibrium.\\

For the chemostat model (\ref{comp}) with two species, one has the following result.
\begin{prop}
\label{prop2species}
Denote by $\lambda_{i}$ the solution (when it exists) of $\mu_{i}(\lambda_{i})=D$.
Under the condition that $S_{in}>\max(\lambda_1,\lambda_2)$, the system (\ref{comp}) admits three equilibrium points given by $E_0 (S_{in},0,0),E_1 (\lambda_1,k_1(S_{in}-\lambda_1),0)$ and $E_2 (\lambda_2,0,k_2(S_{in}-\lambda_2))$. Furthermore, one has
\begin{itemize}
\item when $\lambda_i<\lambda_j$, $E_i$ is locally exponentially stable and $E_0$ and $E_j$ are both unstable.
\item when $\lambda_1=\lambda_2$ then $E_1=E_2$ is a non-hyperbolic equilibrium point.
\end{itemize}
\end{prop}

\textit{Proof.} As before, the equilibrium points are given by the following system
\begin{equation}
\left\{\begin{array}{lll}
\ds -\frac{\mu_1(S)}{k_1}B_1-\frac{\mu_2(S)}{k_2}B_2+D(S_{in}-S)&=0&\vspace{5 pt} \cr
(\mu_1(S)-D)B_1&=0&\vspace{5 pt} \cr
(\mu_2(S)-D)B_2&=0&\vspace{5 pt} \cr 
\end{array}\right. 
\end{equation}
One can find that there exist at most three equilibrium points $E_0(S_{in},0,0), 
E_1(\lambda_1, k_1(S_{in}-\lambda_1),0)$ and $E_2 (\lambda_2, 0, k_2 (S_{in}-\lambda_2))$. \\

For the study of their stability, we write for convenience the dynamics in variables $(Z,B_1,B_2)$ with
$$Z=\frac{B_1}{K_1} +\frac{B_2}{K_2} +S \ . $$
\begin{equation*}\label{compeq}
\left\{\begin{array}{lll}
\ds  \frac{dZ}{dt}&=&D(S_{in}-Z)\vspace{5 pt} \cr
\ds \frac{dB_1}{dt}&=&(\mu_1(Z-\frac{B_1}{K_1} -\frac{B_2}{K_2})-D)B_1\vspace{5 pt} \cr
\ds \frac{dB_2}{dt}&=&(\mu_2(Z-\frac{B_1}{K_1} -\frac{B_2}{K_2})-D)B_2\vspace{5 pt} \cr 
\end{array}\right. 
\end{equation*}
In these coordinates, the Jacobian matrix takes the following form:

$$J(Z,B_1,B_2)=\begin{pmatrix} -D & 0 &0 &\vspace{7pt} \\
                          \star  & -\frac{\mu'_1(S)}{k_1}B_1+\mu_1 (S)-D &                                 -\frac{\mu'_1(S)}{k_2}B_1 &\vspace{7pt} \\
                          
                           \star & -\frac{\mu'_2(S)}{k_1}B_2 &                                          -\frac{\mu'_2(S)}{k_2}B_1+\mu_2 (S)-D&
         \end{pmatrix} \ . $$
At $E_0$, we can check that $J(E_0)$ admits three eigenvalues:
$$v_1(E_{0})=-D<0,\quad v_2(E_{0})=\mu_1(S_{in})-D>0, \quad v_3(E_{0})=\mu_2(S_{in})-D>0 \ . $$
Thus $E_0$ is unstable.
At $E_1$, the eigenvalues are 
$$v_1(E_{1})=-D<0,\quad v_2(E_{1})=-\frac{\mu'_1(\lambda_1)}{k_1}k_1(S_{in}-\lambda_1)<0, \quad  v_3(E_{1})=\mu_2(\lambda_1)-D \ , $$
and symmetrically for $E_2$:
$$v_1(E_{2})=-D<0,\quad v_2(E_{2})=-\frac{\mu'_2(\lambda_2)}{k_2}k_2(S_{in}-\lambda_2)<0, \quad  v_3(E_{2})=\mu_1(\lambda_2)-D \ . $$

One can notice that $\lambda_i<\lambda_j\Leftrightarrow v_3(E_{i)}<0$
and conclude that when $\lambda_i<\lambda_j$, $E_i$ is locally exponentially stable and $E_j$ is unstable. 

For the particular case of $\lambda_1=\lambda_2$ we find that $v_3(E_{1})=v_3(E_{2})=0$. Therefore in this case $E_1$ (that coincides with $E_2$) is a non-hyperbolic equilibrium point.

\end{document}